\input amstex
\documentstyle{amsppt}
\topmatter
\title
The Birch-Swinnerton-Dyer Conjecture
\endtitle
\rightheadtext{The Birch-Swinnerton-Dyer Conjecture}
\author   Jae-Hyun Yang
\endauthor
\abstract {We give a brief description of the Birch-Swinnerton-Dyer
conjecture which is one of the seven Clay problems.}
\endabstract
\magnification =\magstep 1 \pagewidth{12.5cm} \pageheight{17.77cm}
\baselineskip =7mm
\endtopmatter
\document
\NoBlackBoxes

\define\a{\alpha}

\define\G{\Gamma}

\define\la{\lambda}

\define\lrt{\longrightarrow}

\define\psb{\par\smallpagebreak}

\define\pbb{\par\bigpagebreak}

\define\BZ{\Bbb Z}
\define\BC{\Bbb C}
\define\BR{\Bbb R}
\define\BQ{\Bbb Q}
\define\BH{\Bbb H}

\head  \ {\bf 1.\ Introduction} 
\endhead       
\vskip 0.35cm \ \ \ On May 24, 2000, the Clay Mathematics
Institute (CMI for short) announced that it would award prizes of
1 million dollars each for solutions to seven mathematics
problems. These seven problems are \par\medpagebreak Problem 1.
The ``P versus NP" Problem\,:\psb Problem 2. The Riemann
Hypothesis\,:\psb Problem 3. The Poincar{\'e} Conjecture\,: \psb
Problem 4. The Hodge Conjecture\,: \psb Problem 5. The
Birch-Swinnerton-Dyer Conjecture\,:\psb Problem 6. The
Navier-Stokes Equations\,: Prove or disprove the existence and
\psb \qquad\quad\qquad\ smoothness of solutions to the three
dimensional Navier-Stokes \psb\qquad\quad\qquad\ equations. \psb
Problem 7. Yang-Mills Theory\,: Prove that quantum Yang-Mills
fields exist \psb\qquad\qquad\ \ \ \  and have a mass gap. \pbb
Problem 1 is arisen from theoretical computer science, Problem 2
and Problem 5 from number theory, Problem 3 from topology, Problem
4 from algebraic geometry and topology, and finally problem 6 and
7 are related to physics. For more details on some stories about
these problems, we refer to Notices of AMS, vol. 47, no. 8,
pp.\,877-879\,(September 2000) and the homepage of CMI.
\par\medpagebreak In this paper, I will explain Problem 5, that is, the
Birch-Swinnerton-Dyer conjecture which was proposed by the English
mathematicians, B. Birch and H. P. F. Swinnerton-Dyer around 1960
in some detail. This conjecture says that if $E$ is an elliptic
curve defined over $\BQ$, then the algebraic rank of $E$ equals
the analytic rank of $E.$ Recently the Taniyama-Shimura conjecture
stating that any elliptic curve defined over $\BQ$ is modular was
shown to be true by Breuil, Conrad, Diamond and Taylor\,[BCDT].
This fact shed some lights on the solution of the BSD conjecture.
In the final section, we describe the connection between the
heights of Heegner points on modular curves $X_0(N)$ and Fourier
coefficients of modular forms of half integral weight or of the
Jacobi forms corresponding to them by the Skoruppa-Zagier
correspondence. We would like to mention that we added the nicely
written expository paper [W] of Andrew Wiles about the
Birch-Swinnerton-Dyer Conjecture to the list of the references.
\par\medpagebreak\noindent
{\bf Notations\,:} We denote by $\BQ,\, \BR$ and $\BC$ the fields
of rational numbers, real numbers and complex numbers
respectively. $\BZ$ and $\BZ^+$ denotes the ring of integers and
the set of positive integers respectively.

\vskip 1cm
\head  \ {\bf 2.\ The Mordell-Weil Group} 
\endhead       
\vskip 0.35cm \ \ \ A curve $E$ is said to be an {\it elliptic
curve} over $\BQ$ if it is a nonsingular projective curve of genus
1 with its affine model $$y^2 =f(x),\tag2.1$$ where $f(x)$ is a
polynomial of degree 3 with integer coefficients and with 3
distinct roots over $\BC$. An elliptic curve over $\BQ$ has an
abelian group structure with distinguished element $\infty$ as an
identity element. The set $E(\BQ)$ of rational points given by
$$E(\BQ)=\left\{\,(x,y)\in\BQ^2\,|\ y^2=f(x)\,\right\}\cup
\{\infty\} \tag2.2$$ also has an abelian group
structure.\par\medpagebreak L. J. Mordell\,(1888-1972) proved the
following theorem in 1922.\psb\noindent {\bf Theorem A\,(Mordell,
1922)}. $E(\BQ)$ is finitely generated, that is, $$E(\BQ)\cong
\BZ^r \oplus E_{\text{tor}}(\BQ),$$ where $r$ is a nonnegative
integer and $E_{\text{tor}}(\BQ)$ is the torsion subgroup of
$E(\BQ)$.
\define\pmm{\par\medpagebreak}
\pmm\noindent {\bf Definition 1.} Around 1930, A.
Weil\,(1906-1998) proved the set $A(\BQ)$ of rational points on an
abelian variety $A$ defined over $\BQ$ is finitely generated. An
elliptic curve is an abelian variety of dimension one. Therefore
$E(\BQ)$ is called the {\it Mordell}-{\it Weil group} and the
integer $r$ is said to be the {\it algebraic rank} of $E$. \pmm In
1977, B. Mazur\,(1937-\ )[Ma1] discovered the structure of the
torsion subgroup $E_{\text{tor}}(\BQ)$ completely using a deep
theory of elliptic modular curves. \psb\noindent {\bf Theorem
B\,(Mazur, 1977).} Let $E$ be an elliptic curve defined over
$\BQ$. Then the torsion subgroup $E_{\text{tor}}(\BQ)$ is
isomorphic to the following 15 groups $$\BZ/n\BZ\quad (1\leq n\leq
10,\ n=12),$$ $$\BZ/2\BZ \times \BZ/2n\BZ\quad (1\leq n\leq 4).$$
\indent E. Lutz\,(1914-?) and T. Nagell\,(1895-?) obtained the
following result independently. \pmm\noindent {\bf Theorem
C\,(Lutz, 1937;\ Nagell, 1935).}  Let $E$ be an elliptic curve
defined over $\BQ$ given by $$E\,:\quad y^2=x^2+ax+b,\quad
a,b\in\BZ,\ 4a^3+27b^2\not= 0.$$ Suppose that $P=(x_0,y_0)$ is an
element of the torsion subgroup $E_{\text{tor}}(\BQ)$. Then \psb
(a) $\quad x_0,y_0\in \BZ$, and \psb (b) $2P=0$\ \ or\ \ $y_0^2 |
(4a^3+27b^2).$ \pmm We observe that the above theorem gives an
effective method for bounding $E_{\text{tor}}(\BQ)$. According to
Theorem B and C, we know the torsion part of $E(\BQ)$
satisfactorily. But we have no idea of the free part of $E(\BQ)$
so far. As for the algebraic rank $r$ of an elliptic curve $E$
over $\BQ$, it is known by J.-F. Mestre in 1984 that values as
large as 14 occur. Indeed, the elliptic curve defined by
$$y^2=x^3-35971713708112\,x+85086213848298394000$$ has its
algebraic rank 14. \pmm\noindent {\bf Conjecture D.} Given a
nonnegative integer $n$, there is an elliptic curve $E$ over $\BQ$
with its algebraic rank $n$. \pmm The algebraic rank of an
elliptic curve is an invariant under the isogeny. Here an isogeny
of an elliptic curve $E$ means a holomorphic map $\varphi:
E(\BC)\lrt E(\BC)$ satisfying the condition $\varphi(0)=0.$

\vskip 1cm
\head  \ {\bf 3.\ Modular Elliptic Curves} 
\endhead       
\vskip 0.35cm \ \ \ For a positive integer $N\in\BZ^+,$ we let
$$\Gamma_0(N):=\left\{ \pmatrix a & b\\ c & d \endpmatrix \in
SL(2,\BZ)\,\big| \ N|c\ \right\}$$ be the Hecke subgroup of
$SL(2,\BZ)$ of level $N$. Let $\BH$ be the upper half plane. Then
$$Y_0(N)=\BH/\Gamma_0(N)$$ is a noncompact surface, and
$$X_0(N)=\BH\cup \BQ\cup \{\infty\} /\Gamma_0(N)\tag3.1$$ is a
compactification of $Y_0(N).$ We recall that a {\it cusp form} of
weight $k\geq 1$ and level $N\geq 1$ is a holomorphic function $f$
on $\BH$ such that for all $\pmatrix a & b\\ c & d \endpmatrix \in
\Gamma_0(N)$ and for all $z\in\BH$, we have
$$f((az+b)/(cz+d))=(cz+d)^k f(z)$$ and $|f(z)|^2 (\text{Im}\,z)^k$
is bounded on $\BH$. We denote the space of all cusp forms of
weight $k$ and level $N$ by $S_k(N)$. If $f\in S_k(N)$, then it
has a Fourier expansion $$f(z)=\sum_{n=1}^{\infty}c_n(f) q^n,\quad
q:=e^{2\pi i z}$$ convergent for all $z\in\BH.$ We note that there
is no constant term due to the boundedness condition on $f$. Now
we define the $L$-series $L(f,s)$ of $f$ to be
$$L(f,s)=\sum_{n=1}^{\infty}c_n(f)\,n^{-s}.\tag3.2$$ \indent For
each prime $p\not| N$, there is a linear operator $T_p$ on
$S_k(N)$, called the Hecke operator, defined by
$$(f|{T_p})(z)=p^{-1}\sum_{i=0}^{p-1}f((z+i)/p)+p^{k-1}(cpz+d)^k\cdot
f((apz+d)/(cpz+d))$$ for any $\pmatrix a & b\\ c & d \endpmatrix
\in SL(2,\BZ)$ with $c\equiv 0\, (N)$ and $d\equiv p\,(N).$ The
Hecke operators $T_p$ for $p\not| N$ can be diagonalized on the
space $S_k(N)$ and a simultaneous eigenvector is called an {\it
eigenform.} If $f\in S_k(N)$ is an eigenform, then the
corresponding eigenvalues, $a_p(f)$, are algebraic integers and we
have $c_p(f)=a_p(f)\,c_1(f).$\psb Let $\la$ be a place of the
algebraic closure ${\bar {\BQ}}$ in $\BC$ above a rational prime
$\ell$ and ${\bar {\BQ}}_{\la}$ denote the algebraic closure of
$\BQ_{\ell}$ considered as a ${\bar {\BQ}}$-algebra via $\la$. It
is known that if $f\in S_k(N),$ there is a unique continuous
irreducible representation $$\rho_{f,\la}:\text{Gal}({\bar
{\BQ}}/\BQ)\lrt GL_2({\bar{\BQ}}_{\la})\tag3.3$$ such that for any
prime ${p\not| N}  {\ell}$,\ $\rho_{f,\la}$ is unramified at $p$
and $\text{tr}\,\rho_{f,\la}(\text{Frob}_p)=a_p(f).$ The existence
of $\rho_{f,\la}$ is due to G. Shimura\,(1930- ) if $k=2$ [Sh], to
P. Deligne\,(1944- ) if $k>2$ [D] and to P. Deligne and J.-P.
Serre\,(1926- ) if $k=1$ [DS]. Its irreducibility is due to Ribet
if $k>1$ [R], and to Deligne and Serre if $k=1$ [DS]. Moreover
$\rho_{f,\la}$ is odd and potentially semi-stable at $\ell$ in the
sense of Fontaine. We may choose a conjugate of $\rho_{f,\la}$
which is valued in $GL_2({\Cal O}_{{\bar{\BQ}}_{\la}})$, and
reducing modulo the maximal ideal and semi-simplifying yields a
continuous representation $${\bar\rho}_{f,\la}:\text{Gal}({\bar
{\BQ}}/\BQ)\lrt GL_2({\bar{\Bbb F}}_{\ell}),\tag3.4$$ which, up to
isomorphism, does not depend on the choice of conjugate of
$\rho_{f,\la}$. \pmm\noindent {\bf Definition 2.} Let
$\rho:\text{Gal}({\bar {\BQ}}/\BQ)\lrt GL_2({\bar{\BQ}}_{\ell})$
be a continuous representation which is unramified outside
finitely many primes and for which the restriction of $\rho$ to a
decomposition group at $\ell$ is potentially semi-stable in the
sense of Fontaine. We call $\rho$ {\it modular} if $\rho$ is
isomorphic to $\rho_{f,\la}$ for some eigenform $f$ and some $\la
| \ell.$\pmm\noindent {\bf Definition 3.} An elliptic curve $E$
defined over $\BQ$ is said to be {\it modular} if there exists a
surjective holomorphic map $\varphi: X_0(N)\lrt E(\BC)$ for some
positive integer $N$.\psb Recently C. Breuil, B. Conrad, F.
Diamond and R. Taylor [BCDT] proved that the Taniyama-Shimura
conjecture is true. \psb\noindent {\bf Theorem E\,([BCDT], 2001).}
An elliptic curve defined over $\BQ$ is modular. \pmm Let $E$ be
an elliptic curve defined over $\BQ$. For a positive integer
$n\in\BZ^+$, we define the isogeny $[n]:E(\BC)\lrt E(\BC)$ by
$$[n]P:=nP=P+\cdots+P\ (n\ \text{times}),\quad P\in
E(\BC).\tag3.5$$ For a negative integer $n$, we define the isogeny
$[n]:E(\BC)\lrt E(\BC)$ by $[n]P:=-[-n]P,\ P\in E(\BC)$, where
$-[-n]P$ denotes the inverse of the element $[-n]P$. And
$[0]:E(\BC)\lrt E(\BC)$ denotes the zero map. For an integer
$n\in\BZ,\ [n]$ is called the multiplication-by-$n$ homomorphism.
The kernel $E[n]$ of the isogeny $[n]$ is isomorphic to
$\BZ/n\BZ\oplus \BZ/n\BZ.$ Let
$$\text{End}(E)=\left\{\varphi:E(\BC)\lrt E(\BC),\ \text{an\
isogeny}\,\right\}$$ be the endomorphism group of $E$. An elliptic
curve $E$ over $\BQ$ is said to have {\it complex multiplication}
(or CM for short) if $$\text{End}(E)\not\subseteq \BZ\cong \left\{
[n]\vert\ n\in\BZ\,\right\},$$ that is, there is a nontrivial
isogeny $\varphi:E(\BC)\lrt E(\BC)$ such that $\varphi\not= [n]$
for all integers $n\in\BZ.$ Such an elliptic curve is called a CM
{\it curve}. For most of elliptic curves $E$ over $\BQ$, we have
$\text{End}(E)\cong \BZ.$

\vskip 1cm
\head  \ {\bf 4.\ The $L$-Series of an Elliptic Curve} 
\endhead       
\vskip 0.35cm \ \ \ Let $E$ be an elliptic curve over $\BQ$. The
$L$-series $L(E,s)$ of $E$ is defined  as the product of the local
$L$-factors\,: $$L(E,s)=\prod_{p| \Delta_E} (1-a_p p^{-s})^{-1}
\cdot \prod_{{p \not|} \Delta_E} (1-a_p
p^{-s}+p^{1-2s})^{-1},\tag4.1$$ where $\Delta_E$ is the
discriminant of $E$, $p$ is a prime, and if $p \not| \Delta_E$,
$$a_p:=p+1-|{\bar E}({\Bbb F}_p)|,$$ and if $p| \Delta_E,$ we
set $a_p:=0,\,1,\,-1$ if the reduced curve ${\bar E}/{\Bbb F}_p$
has a cusp at $p$, a split node at $p$, and a nonsplit node at $p$
respectively. Then $L(E,s)$ converges absolutely for
$\text{Re}\,s>{\frac 32}$ from the classical result that $|a_p|< 2
\sqrt{p}$ for each prime $p$ due to H. Hasse\,(1898-1971) and is
given by an absolutely convergent Dirichlet series. We remark that
$x^2-a_p x+p$ is the characteristic polynomial of the Frobenius
map acting on ${\bar E}({\Bbb F}_p)$ by $(x,y)\mapsto (x^p,y^p).$
\pmm\noindent {\bf Conjecture F.} Let $N(E)$ be the conductor of
an elliptic curve $E$ over $\BQ$\,([S],\,p.\,361). We set
$$\Lambda(E,s):=N(E)^{s/2}\,(2\pi)^{-s}\,\Gamma(s)\, L(E,s),\quad
\text{Re}\,s>{\frac 32}.$$ Then $\Lambda(E,s)$ has an analytic
continuation to the whole complex plane and satisfies the
functional equation $$\Lambda(E,s)=\epsilon\, \Lambda(E,2-s),\quad
\epsilon=\pm 1.$$ \indent The above conjecture is now true because
the Taniyama-Shimura conjecture is true\,(cf. Theorem E). We have
some knowledge about analytic properties of $E$ by investigating
the $L$-series $L(E,s)$. The order of $L(E,s)$ at $s=1$ is called
the {\it analytic rank} of $E$. \psb Now we explain the connection
between the modularity of an elliptic curve $E$, the modularity of
the Galois representation and the $L$-series of $E$. For a prime
$\ell$, we let $\rho_{E,\ell}$\,(resp. ${\bar\rho}_{E,\ell})$
denote the representation of $\text{Gal}({\bar\BQ}/\BQ)$ on the
$\ell$-adic Tate module\,(resp.\ the $\ell$-torsion) of
$E({\bar\BQ}).$ Let $N(E)$ be the conductor of $E$. Then it is
known that the following conditions are equivalent\,:\psb (1) The
$L$-function $L(E,s)$ of $E$ equals the $L$-function $L(f,s)$ for
some eigenform \psb \qquad$f$. \psb (2) The $L$-function $L(E,s)$
of $E$ equals the $L$-function $L(f,s)$ for some eigenform
\psb\qquad $f$ of weight 2 and level $N(E)$.\psb (3) For some
prime $\ell$, the representation $\rho_{E,\ell}$ is modular.\psb
(4) For all primes $\ell$, the representation $\rho_{E,\ell}$ is
modular.\psb (5) There is a non-constant holomorphic map
$X_0(N)\lrt E(\BC)$ for some positive \psb\qquad integer $N$.\psb
(6) There is a non-constant morphism $X_0(N(E))\lrt E$ which is
defined over \psb\qquad $\BQ.$\psb (7) $E$ admits a hyperbolic
uniformization of arithmetic type\,(cf.\,[Ma2] and [Y1]).

\vskip 1cm
\head {\bf 5.\ The Birch-Swinnerton-Dyer conjecture} 
\endhead       
\vskip 0.35cm \ \ \ Now we state the BSD conjecture.\psb\noindent
{\bf The BSD Conjecture.} Let $E$ be an elliptic curve over $\BQ$.
Then the algebraic rank of $E$ equals the analytic rank of $E.$
\pmm I will describe some historical backgrounds about the BSD
conjecture. Around 1960, Birch\,(1931- ) and
Swinnerton-Dyer\,(1927- ) formulated a conjecture which determines
the algebraic rank $r$ of an elliptic curve $E$ over $\BQ$. The
idea is that an elliptic curve with a large value of $r$ has a
large number of rational points and should therefore have a
relatively large number of solutions modulo a prime $p$ on the
average as $p$ varies. For a prime $p$, we let $N(p)$ be the
number of pairs of integers $x,y\,(\text{mod}\,p)$ satisfying
(2.1) as a congruence (mod $p$). Then the BSD conjecture in its
crudest form says that we should have an asymptotic formula
$$\prod_{p<x}{{N(p)+1}\over p}\ \sim C \ \,(\text{log}\,p)^r \quad
\text{as}\ x\lrt\infty\tag5.1$$ for some constant $C>0.$ If the
$L$-series $L(E,s)$ has an analytic continuation to the whole
complex plane\,(this fact is conjectured; cf.\, Conjecture F),
then $L(E,s)$ has a Taylor expansion  $$L(E,s)=c_0
(s-1)^m+c_1(s-1)^{m+1}+\cdots$$ at $s=1$ for some non-negative
integer $m\geq 0$ and constant $c_0\not= 0.$ The BSD conjecture
says that the integer $m$, in other words, the analytic rank of
$E$, should equal the algebraic rank $r$ of $E$ and furthermore
the constant $c_0$ should be given by $$c_0=\lim_{s\rightarrow
1}{{L(E,s)}\over {(s-1)^m}}=\a\cdot R\cdot
|E_{\text{tor}}(\BQ)|^{-1}\cdot \Omega\cdot S,\tag5.2$$ where
$\a>0$ is a certain constant, $R$ is the elliptic regulator of
$E,\ |E_{\text{tor}}(\BQ)|$ denotes the order of the torsion
subgroup $E_{\text{tor}}(\BQ)$ of $E(\BQ)$, $\Omega$ is a simple
rational multiple (depending on the bad primes) of the elliptic
integral $$\int_{\gamma}^{\infty}{{dx}\over {\sqrt{f(x)}}}\qquad
(\gamma=\text{the\ largest\ root\ of}\ f(x)=0)$$ and $S$ is an
integer square which is supposed to be the order of the
Tate-Shafarevich group $\text{III}(E)$ of $E$. \psb The
Tate-Shafarevich group $\text{III}(E)$ of $E$ is a very intersting
subject to be investigated in the future. Unfortunately
$\text{III}(E)$ is still not known to be finite. So far an
elliptic curve whose Tate-Shafarevich group is infinite has not
been discovered. So many mathematicians propose the following.
\psb\noindent {\bf Conjecture G.} The Tate-Shafarevich group
$\text{III}(E)$ of $E$ is finite. \pmm There are some evidences
supporting the BSD conjecture. I will list these evidences
chronologically. \psb\noindent {\bf Result
1\,}(Coates-Wiles\,[CW],\,1977). Let $E$ be a CM curve over $\BQ$.
Suppose that the analytic rank of $E$ is zero. Then the algebraic
rank of $E$ is zero. \psb\noindent {\bf Result
2\,}(Rubin\,[R],\,1981). Let $E$ be a CM curve over $\BQ$. Assume
that the analytic rank of $E$ is zero. Then the Tate-Shafarevich
group $\text{III}(E)$ of $E$ is finite. \psb\noindent {\bf Result
3\,}(Gross-Zagier\,[GZ],\,1986\,;\ [BCDT],\,2001). Let $E$ be an
elliptic curve over $\BQ$. Assume that the analytic rank of $E$ is
equal to one and $\epsilon=-1$\,(cf.\,Conjecture F). Then the
algebraic rank of $E$ is equal to or bigger than one.
\psb\noindent {\bf Result 4\,}(Gross-Zagier\,[GZ],\,1986). There
exists an elliptic curve $E$ over $\BQ$ such that
$\text{rank}\,E(\BQ)=\text{ord}_{s=1} L(E,s)=3$. For instance, the
elliptic curve ${\tilde E}$ given by $${\tilde E}\ :\quad
-139\,y^2=x^3+10\,x^2-20\,x+8$$ satisfies the above property.
\psb\noindent {\bf Result 5\,}(Kolyvagin\,[K],\,
1990\,:\,Gross-Zagier\,[GZ],\,1986\,:\,Bump-Friedberg-Hoffstein\,\psb\noindent
[BFH],
\,1990\,:\,Murty-Murty\,[MM],\,1990\,:\,[BCDT],\,2001). Let $E$ be
an elliptic curve over $\BQ$. Assume that the analytic rank of $E$
is 1 and $\epsilon=-1$. Then algebraic rank of $E$ is equal to 1.
\psb\noindent {\bf Result 6\,}(Kolyvagin\,[K],\,
1990\,:\,Gross-Zagier\,[GZ],\,1986\,:\,Bump-Friedberg-Hoffstein\,\psb\noindent[BFH],
\,1990\,:\,Murty-Murty\,[MM],\,1990\,:\,[BCDT],\,2001). Let $E$ be
an elliptic curve over $\BQ$. Assume that the analytic rank of $E$
is zero and $\epsilon=1$. Then algebraic rank of $E$ is equal to
zero. \pmm Cassels proved the fact that if an elliptic curve over
$\BQ$ is isogeneous to another elliptic curve $E'$ over $\BQ$,
then the BSD conjecture holds for $E$ if and only if th e BSD
conjecture holds for $E'$.

\vskip 1cm
\head {\bf 6.\ Jacobi Forms and Heegner Points} 
\endhead       
\vskip 0.35cm \indent In this section, I shall describe the result
of Gross-Kohnen-Zagier\,[GKZ] roughly.\psb First we begin with
giving the definition of Jacobi forms. By definition a Jacobi form
of weight $k$ and index $m$ is a holomorphic complex valued
function $\phi(z,w)\,(z\in \BH,\,z\in\BC)$ satisfying the
transformation formula $$ \align \phi\left( {{az+b}\over {cz+d}},
{{w+\la z+\mu}\over{cz+d}}\right)=&e^{-2\pi i\left\{ cm(w+\la
z+\mu)^2 (cz+d)^{-1}-m(\la^2 z+2\la w)\right\}} \\ &\ \times
(cz+d)^k \,\phi(z,w) \tag6.1\endalign$$ for all $\pmatrix a & b\\
c &d\endpmatrix\in SL(2,\BZ)$ and $(\la,\mu)\in\BZ^2$ having a
Fourier expansion of the form $$\phi(z,w)=\sum_{\scriptstyle
n,r\in \BZ^2 \atop\scriptstyle r^2\leq 4mn} c(n,r)\,e^{2\pi
i(nz+rw)}.\tag6.2$$ We remark that the Fourier coefficients
$c(n,r)$ depend only on the discrimnant $D=r^2-4mn$ and the
residue $r\,(\text{mod}\ 2m).$ From now on, we put
$\Gamma_1:=SL(2,\BZ).$ We denote by $J_{k,m}(\Gamma_1)$ the space
of all Jacobi forms of weight $k$ and index $m$. It is known that
one has the following isomorphisms $$[\Gamma_2,k]^M\cong
J_{k,1}(\G_1)\cong M_{k-{\frac 12}}^+(\G_0(4))\cong
[\G_1,2k-2],\tag6.3$$ where $\G_2$ denotes the Siegel modular
group of degree 2, $[\G_2,k]^M$ denotes the Maass space introduced
by H. Maass\,(1911-1993)\,(cf.\,[M1-3]), $M_{k-{\frac
12}}^+(\G_0(4))$ denotes the Kohnen space introduced by W.
Kohnen\,[Koh] and $[\G_1,2k-2]$ denotes the space of modular forms
of weight $2k-2$ with respect to $\G_1$. We refer to [Y1] and
[Y3],\,pp.\,65-70 for a brief detail. The above isomorphisms are
compatible with the action of the Hecke operators. Moreover,
according to the work of Skoruppa and Zagier [SZ], there is a
Hecke-equivariant correspondence between Jacobi forms of weight
$k$ and index $m$, and certain usual modular forms of weight
$2k-2$ on $\G_0(N).$ \psb Now we give the definition of Heegner
points of an elliptic curve $E$ over $\BQ$. By [BCDT], $E$ is
modular and hence one has a surjective holomorphic map
$\phi_E:X_0(N)\lrt E(\BC).$ Let $K$ be an imaginary quadratic
field of discriminant $D$ such that every prime divisor $p$ of $N$
is split in $K$. Then it is easy to see that $(D,N)=1$ and $D$ is
congruent to a square $r^2$ modulo $4N$. Let $\Theta$ be the set
of all $z\in\BH$ satisfying the following conditions
$$az^2+bz+c=0,\quad a,b,c\in\BZ,\ N|a,$$ $$b\equiv
r\,\,(\text{mod}\,2N),\qquad D=b^2-4ac.$$ Then $\Theta$ is
invariant under the action of $\G_0(N)$ and $\Theta$ has only a
$h_K$ $\G_0(N)$-orbits, where $h_K$ is the class number of $K$.
Let $z_1,\cdots,z_{h_K}$ be the representatives for these
$\G_0(N)$-orbits. Then $\phi_E(z_1),\cdots,\phi_E(z_{h_K})$ are
defined over the Hilbert class field $H(K)$ of $K$, i.e., the
maximal everywhere unramified extension of $K$. We define the
Heegner point $P_{D,r}$ of $E$ by
$$P_{D,r}=\sum_{i=1}^{h_K}\phi_E(z_i).\tag6.4$$ We observe that
$\epsilon=1$, then $P_{D,r}\in E(\BQ)$.\psb Let $E^{(D)}$ be the
elliptic curve (twisted from $E$) given by $$E^{(D)}\ :\ \
Dy^2=f(x).\tag6.5$$ Then one knows that the $L$-series of $E$ over
$K$ is equal to $L(E,s)\,L(E^{(D)},s)$ and that $L(E^{(D)},s)$ is
the twist of $L(E,s)$ by the quadratic character of $K/\BQ$.
\pmm\noindent {\bf Theorem
H}\,(Gross-Zagier\,[GZ],\,1986\,;\,[BCDT],\,2001). Let $E$ be an
elliptic curve of conductor $N$ such that $\epsilon =-1.$ Assume
that $D$ is odd. Then
$$L'(E,1)\,L(E^{(D)},1)=c_E\,u^{-2}\,|D|^{-{\frac 12}}\,{\hat
h}(P_{D,r}),\tag6.6$$ where $c_E$ is a positive constant not
depending on $D$ and $r,\ u$ is a half of the number of units of
$K$ and ${\hat h}$ denotes the canonical height of $E$. \psb Since
$E$ is modular by [BCDT], there is a cusp form of weight 2 with
respect to $\G_0(N)$ such that $L(f,s)=L(E,s).$ Let $\phi(z,w)$ be
the Jacobi form of weight 2 and index $N$ which corresponds to $f$
via the Skoruppa-Zagier correspondence. Then $\phi(z,w)$ has a
Fourier series of the form (6.2). \psb B. Gross, W. Kohnen and D.
Zagier obtained the following result. \pmm\noindent {\bf Theorem
I\,}(Gross-Kohnen-Zagier,\,[GKZ]\,;\,{BCDT],\,2001). Let $E$ be a
modular elliptic curve with conductor $N$ and suppose that
$\epsilon=-1,\ r=1.$ Suppose that $(D_1,D_2)=1$ and $D_i\equiv
r_i^2\,(\text{mod}\,4N)\,(i=1,2).$ Then
$$L'(E,1)\,c((r_1^2-D_1)/(4N),r_1)\,c((r_2^2-D_2)/(4N),r_2)\,=\,c_E'<P_{D_1,r_1},
P_{D_2,r_2}>,$$ where $c_E'$ is a positive constant not depending
on $D_1,\,r_1$ and $D_2,\,r_2$ and $<\ ,\ >$ is the height pairing
induced from the N{\' e}ron-Tate height function ${\hat h}$, that
is,  ${\hat h}(P_{D,r})=$\psb\noindent $<P_{D,r},P_{D,r}>$. \psb We see from the
above theorem that the value of $<P_{D_1,r_1}, P_{D_2,r_2}>$ of
two distinct Heegner points is related to the product of the
Fourier coefficients
$c((r_1^2-D_1)/(4N),r_1)\,c((r_2^2-D_2)/(4N),r_2)$ of the Jacobi
forms of weight 2 and index $N$ corresponded to the eigenform $f$
of weight 2 associated to an elliptic curve $E$. We refer to [Y4]
and [Z] for more details.
\pmm\noindent
{\bf Corollary.} There is a point $P_0\in E(\BQ)\otimes_{\BZ}\BR$
such that $$P_{D,r}=c((r^2-D)/(4N),r)P_0$$ for all $D$ and
$r\,(D\equiv r^2\,(\text{mod}\,4N))$ with $(D,2N)=1.$ \psb The
corollary is obtained by combining Theorem H and Theorem I with
the Cauchy-Schwarz inequality in the case of equality.
\pmm\noindent {\bf Remark 4.} R. Borcherds \,[B] generalized the
Gross-Kohnen-Zagier theorem to certain more general quotients of
Hermitian symmetric spaces of high dimension, namely to quotients
of the space associated to an orthogonal group of signature $(2,b)$
by the unit group of a lattice. Indeed he relates the Heegner
divisors on the given quotient space to the Fourier coefficients
of vector-valued holomorphic modular forms of weight $1+{\frac
b2}$.

\vskip 1cm \noindent{\tt Department of Mathematics\par\noindent
Inha University\par\noindent Incheon 402-751\par\noindent Republic
of Korea}
\par\bigpagebreak\noindent \vskip 0.5cm \noindent {\tt email\,:\
jhyang\@inha.ac.kr}

\end{document}